\documentclass[12pt]{amsart}
\usepackage{amssymb, graphicx}
\usepackage[latin1]{inputenc}
\usepackage{amsfonts}
\usepackage{latexsym}
\usepackage{amscd}
\usepackage{enumerate}
\newtheorem{prop}{Proposition}[section]
\newtheorem{lem}[prop]{Lemma}
\newtheorem{thm}[prop]{Theorem}
\newtheorem{cor}[prop]{Corollary}
\newtheorem{definition}[prop]{Definition}
\newtheorem{definitions}[prop]{Definitions}

\theoremstyle{definition}

\newtheorem{rem}[prop]{Remark}

\numberwithin{equation}{section}

\newcommand{\qa}{algebra of quotients }

\newcommand{\wqa}{weak algebra of quotients }

\newcommand{\lia}{left ideally absorbed into }
\newcommand{\ria}{right ideally absorbed into }
\newcommand{\ad}{\mathrm{ad}_}

\newcommand{\idl}{\mathrm{id}^l_}
\newcommand{\idr}{\mathrm{id}^r_}
\newcommand{\id}{\mathrm{id}_}
\newcommand{\lan}{\mathrm{l.ann}_}
\newcommand{\ran}{\mathrm{r.ann}_}
\newcommand{\ann}{\mathrm{Ann}_}

\begin{document}

\title{Associative and Lie algebras of quotients}

\author{Francesc Perera \& Mercedes Siles Molina}
\address{Departament de Matem\`atiques,
Universitat Aut\`onoma de Barcelona, 08193 Be\-lla\-terra, Barcelona,
Spain} \email{perera@mat.uab.es}
\address{Departamento de \'Algebra, Geometr\'\i a y Topolog\'\i a,
Universidad de M\'alaga, 29071 M\'alaga, Spain}
\email{mercedes@agt.cie.uma.es}

\thanks{Research supported by the MCYT and European Regional
Development Fund through Projects BFM2002-01390 and
MTM2004-06580-C02-02, the Comissionat per Universitats i Recerca
de la Generalitat de Catalunya, and the ``Plan Andaluz de
Investigaci\'on y Desarrollo Tecnol\'ogico'', FQM 336.}
\date{\today}

\begin{abstract}
In this paper we examine how the notion of algebra of quotients
for Lie algebras ties up with the corresponding well-known concept
in the associative case. Specifically, we completely characterize
when a Lie algebra $Q$ is an algebra of quotients of a Lie algebra
$L$ in terms of the associative algebras generated by the adjoint
operators of $L$ and $Q$ respectively. In a converse direction, we
also provide with new examples of algebras of quotients of Lie
algebras and these come from associative algebras of quotients. In
the course of our analysis, we make use of the notions of density
and multiplicative semiprimeness to link our results with the
maximal symmetric ring of quotients.
\end{abstract}

\maketitle

\section*{Introduction}

In recent years, there has been a trend to extend notions and
results of algebras of quotients of associative algebras to the
non-associative setting. This has been achieved by a number of
authors, see
e.g.~\cite{mart},~\cite{msm},~\cite{gago},~\cite{monpa},~\cite{agago}.

One of the leitmotifs for carrying out this process is the fact
that, in the associative case, rings of quotients allow a deeper
understanding of certain classes of rings. Thus it is to be
expected that a similar role will be played by their
non-associative siblings.

In this paper we shall be concerned with Lie algebras and their
algebras of quotients. These were introduced by the second author
in~\cite{msm}, following the original pattern of Utumi~\cite{utu} (see below
for the precise definitions). This approach prompts the question
of whether there is a relationship between the associative and Lie
algebras of quotients, beyond the formal analogy of the
definitions. There are at least two ways to analyse this. On the
one hand, any Lie algebra $L$ gives rise to an associative algebra
$A(L)$, which is generated by the adjoint operators given by the
bracket in $L$. If $L\subseteq Q$ are Lie algebras such that $Q$
is an algebra of quotients of $L$, it is then natural to ask
whether $A(Q)$ is, in some sense, an algebra of quotients of
$A(L)$. Of course, one has to circumvent the fact that $A(L)$
might not be well related to $A(Q)$ (it might not be even a
subalgebra of $A(Q)$). In order to deal with this it is natural to
consider the subalgebra $A_Q(L)$ of $A(Q)$ generated by the
adjoint operators coming from elements in $L$, or the subalgebra
$A_0$ of $A(Q)$ whose elements map $L$ into $L$. This contains,
and may not coincide with, $A_Q(L)$.

We prove that, if $L$ and $Q$ are Lie algebras such that $Q$ is an
algebra of quotients of $L$, then $A(Q)$ is a left quotient
algebra of $A_0$. This is one of the main themes of Section 2 and is accomplished in
Corollary~\ref{leftquotient}. The exact relationship between the
property of $Q$ being a quotient algebra of $L$ and the
corresponding property in terms of $A(Q)$ is of a more technical
nature and is established in Theorem~\ref{elteorema}. This result
is also important for the study of dense extensions of Lie
algebras (see below).

Starting from the other endpoint, we may consider an extension of
(semiprime) associative algebras $A\subseteq Q$ such that $Q$ is a
quotient algebra of $A$ and $Q\subseteq Q_s(A)$ (the Martindale
symmetric ring of quotients). By considering the corresponding
natural Lie structures on both $A$ and $Q$, it is natural to ask
whether $Q^{(-)}$ is a Lie algebra of quotients of $A^{(-)}$. This
was proved to be the case by the second author
in~\cite[Proposition 2.16]{msm} if moreover $Z(Q)=0$. The general
situation is much more intricate and thus requires a deeper
analysis. Our main result in this direction uses arguments that go back to Herstein and asserts that
$Q^{(-)}/Z(Q^{(-)})$ is a Lie quotient algebra of $A^{(-)}/Z(A^{(-)})$ whenever $A\subseteq
Q\subseteq Q_s(A)$ (see Theorem~\ref{liecenter}). As a consequence we obtain in
Corollary~\ref{liequocenter} that
$[Q^{(-)},Q^{(-)}]/Z([Q^{(-)},Q^{(-)}])$ is a quotient algebra of
$[A^{(-)},A^{(-)}]/Z([A^{(-)},A^{(-)}])$, under the same
assumptions on $A$ and $Q$.

Going back to the relationship between a Lie algebra $L$ and the
associative algebra $A(L)$, we examine in Section 3 the behaviour of this
association under extensions. More concretely, if $L\subseteq Q$
are Lie algebras, we consider when the natural correspondence
$A(L)\to A(Q)$ given by the change of domain is an algebra map.
This amounts to requiring that the extension $L\subseteq Q$ is
dense in the sense of Cabrera (see~\cite{Ca}).

The first examples of dense extensions were given in the context
of multiplicatively semiprime algebras. Recall that an algebra $A$
(associative or not) is said to be multiplicatively semiprime if
$A$ and the multiplication algebra $M(A)$ (generated by the right
and left multiplication operators of $A$ together with the
identity) are both semiprime (see~\cite{CaM}, \cite{CabM},
\cite{CaCa}, \cite{CCLM} among others). Roughly, an extension
$A\subseteq B$ of algebras is dense if every non-zero element in
$M(B)$ remains non-zero when restricted to $A$. Cabrera proved
in~\cite{Ca} that every essential ideal of a multiplicatively
semiprime algebra is dense.

By using our techniques, we are able to find new and significant
instances where dense extensions naturally appear. For example, if
$A\subseteq Q\subseteq
Q_{\mathrm{max}}^l(A)$ are associative algebras, where $Q_{\mathrm{max}}^l(A)$ is the maximal left quotient ring of $A$, then this is a dense extension, and the
corresponding extensions of Lie algebras $A^{(-)}/Z(A)\subseteq
Q^{(-)}/Z(Q)$ and $
[A^{(-)},A^{(-)}]/Z([A^{(-)},A^{(-)}])\subseteq
[Q^{(-)},Q^{(-)}]/Z([Q^{(-)},Q^{(-)}]) $ are also dense (see
Lemma~\ref{denseassoc} and Proposition~\ref{denselie}).

Our considerations finally lead to the conclusion that, if
$L\subseteq Q$ is a dense extension of Lie algebras with $Q$
multiplicatively semiprime, then the associative algebras $A_0$
and $A(Q)$ are indistinguishable under the formation of the
maximal symmetric ring of quotients (in the sense of
Schelter~\cite{sch} and Lanning~\cite{lann}), see
Theorem~\ref{20/02/05,7}.

This paper was done during visits of the first author to the
Universidad de Málaga and of the second author to the Universitat
Autònoma de Barcelona. Both authors wish to thank the respective
host centers for their hospitality.

\section{Preliminaries}

In this paper we will deal with algebras over an arbitrary
unital and commutative ring of scalars $\Phi$. In the case our algebras are associative, they need not be unital. First of all we will give some definitions and basic notation.

A $\Phi$-module $L$ with a bilinear map $[\ ,\ ]\colon L\times L
\to L$, denoted by $(x, y)\mapsto [x, y]$ and called the
\emph{bracket} of $x$ and $y$, is called a  \emph{Lie algebra} over
$\Phi$ if the following axioms are satisfied:
\begin{enumerate}[(i)]\itemsep=2mm
\item $[x, x]=0$,
\item $[x, [y, z]]+[y, [z, x]]+[z, [x, y]]=0$ (\emph{Jacobi identity}),
\end{enumerate}
\par
\noindent for every $x, y, z$ in $L$.

Every associative algebra $A$ gives rise to a Lie algebra
$A^{(-)}$ by considering the same module structure and bracket
given by $[x, y]=xy-yx$.

Any element $x$ of a Lie algebra $L$ determines a map $ \ad
x\colon L \to L$ defined by $\ad x(y)=[x, y]$, which is a
derivation of the Lie algebra $L$. We shall denote by $A(L)$ the
associative subalgebra (possibly without identity) of
$\mathrm{End}(L)$ generated by the elements $\ad x$ for $x$ in
$L$. This has been referred to as the \emph{multiplication ideal}
in the literature (see, e.g.~\cite[p. 3493]{CaCa}).

An element $x$ in a Lie algebra $L$ is an \emph{absolute zero
divisor} if $(\ad x)^2=0$. This is equivalent to saying that $[[L,
x], x]= 0$.
The algebra $L$ is said to be \emph{non-degenerate} (or strongly
nondegenerate according to Kostrikin) if it does not contain
non-zero absolute zero divisors.

Given an ideal $I$ in a Lie algebra $L$, we write $I^2=[I,I]$,
which is again an ideal. We say that $L$ is \emph{semiprime} if we
have $I^2\neq 0$ whenever $I$ is a non-zero ideal. In other words,
$L$ has no abelian ideals. Clearly, non-degenerate Lie algebras
are semiprime, while it is possible to find examples where the
converse implication does not hold.

For a subset $X$ of a Lie algebra $L$ recall that the set
\[
\ann{}(X)=\{ a \in L \ \vert \ [a, x]=0 \text{ for every } x
\text{ in } X\}
\]
is called the \emph{annihilator} of $X$. This will be denoted by
$\ann L(X)$ when it is necessary to emphasize the dependence on
$L$. If $X=L$, then $\ann{}(L)$ is called the \emph{centre} of $L$
and usually denoted by $Z(L)$. In the case that $L=A^{(-)}$ for an
associative algebra $A$, then $Z(A^{(-)})$ agrees with the
associative center $Z(A)$ of $A$. It is easy to check (by using
the Jacobi identity) that $\ann{}(X)$ is an ideal of $L$ when $X$
is an ideal of $L$. Every element of $\ann{}(L)$ will be called a
\emph{total zero divisor}.

By an \emph{extension of Lie algebras} $L\subseteq Q$ we will mean
that $L$ is a (Lie) subalgebra of the Lie algebra $Q$.

Let $L\subseteq Q$ be an extension of Lie algebras and let
$A_Q(L)$ be the associative subalgebra of $A(Q)$ generated by
$\{\ad{x} : x \in L\}$.

Recall that, given an associative algebra $A$ and a subset $X$ of
$A$, we define the \emph{right annihilator of} $X$ \emph{in} $A$
as
\[
\ran A(X)=\{a\in A\mid Xa=0\}\,,
\]
which is always a right ideal of $A$ (and two-sided if $X$ is a
right $A$-module). One similarly defines the \emph{left
annihilator}, which shall be denoted by $\lan A(X)$.
\begin{lem}
\label{26/11/04,3} Let $I$ be an ideal of a Lie algebra $L$ with
$Z(L)=0$. Then $\ann{L}(I)=0$ if and only if
$\ran{A(L)}(A_L(I))=0$.
\end{lem}
\begin{proof}
Suppose first $\ann{L}(I)=0$. Let $\mu\in\ran{A(L)}(A_L(I))$. Then
$\ad{y}\mu=0$ for all $y$ in $I$. In particular, if $x\in L$ we
get $0=\ad y\mu(x)=[y,\mu(x)]$, and this implies that
$\mu(x)\in\ann{L}(I)=0$. Hence $\mu=0$.

Conversely, suppose $\ran{A(L)}(A_L(I))=0$. If $x\in\ann{L}(I)$,
$y\in I,$ $z\in L$, then $\ad{y}\ad{x}z= [y,[x, z]]= [[y, x], z]+
[x, [y, z]]=0$, so  $\ad{x}\in \ran{A(L)}(A_L(I))=0$. Since by
assumption $Z(L)=0$ we obtain $x=0$.
\end{proof}
For a subset $X$ of an associative algebra $A$, denote by
$\idl{A}(X)$, $\idr{A}(X)$ and $\id{A}(X)$ the left, right and two
sided ideal of $A$, respectively, generated by $X$. When it is
clear from the context, the reference to the algebra where these
ideals sit into will be omitted.
\begin{lem}
\label{26/11/04,1} Let $I$ be an ideal of a Lie algebra $L$ and
suppose that $Q$ is an overalgebra of $L$. Then $\idl
{A_Q(L)}(A_Q(I))=\idr {A_Q(L)}(A_Q(I))=\id {A_Q(L)}(A_Q(I))$.
\end{lem}
\begin{proof}
Use induction and notice that given $x$ in $L$ and $y$ in $I$ we
have $\ad x \ad y= \ad{[x, y]}+ \ad y \ad x$.
\end{proof}
\begin{lem}
\label{26/11/04,2} Let $I$ be an ideal of a Lie algebra $L$ and
suppose that $Q$ is an overalgebra of $L$. Write $\id{}(A(I))$ to
denote the ideal of $A_{Q}(L)$ generated by $A_Q(I)$. Then
\begin{enumerate}[{\rm (i)}]\itemsep=2mm
\item $\ran{A_Q(L)}(\id{}(A(I)))= \ran{A_Q(L)}(A_Q(I))$.
\item $\lan{A_Q(L)}(\id{}(A(I)))= \lan{A_Q(L)}(A_Q(I))$.
\end{enumerate}
\end{lem}
\begin{proof} To see (i), it is enough to prove that
$\ran{A_Q(L)}(A_Q(I))\subseteq \ran{A_Q(L)}(\id{}(A(I)))$ because
the converse inclusion is obvious. Let
$\lambda\in\ran{A_Q(L)}(A_Q(I))$. By Lemma~\ref{26/11/04,1} we
know that, if $\mu\in \id{}(A(I))$ there exist a natural number
$n$, elements $x_{1, i}, \ldots, x_{r_i, i}$ in $L$, and $y_{1,
i}, \ldots, y_{s_i, i}$ in $I$ with $0\leq r_i\in\mathbb{N}$ for
all $i$ and $\emptyset \neq \{s_1, \ldots, s_n\}\subseteq
\mathbb{N}$, such that $\mu=\sum_{i=1}^{n}\ad{x_{1, i}}\cdots
\ad{x_{r_i, i}}\ad{y_{1, i}}\cdots \ad{y_{s_i, i}}$. Since
$\ad{y_{s_i, i}}\lambda=0$, we see that $\mu\lambda=0$.

A similar argument establishes (ii).
\end{proof}

\section{Algebras of quotients: The Lie and associative cases}

Inspired by the notion of ring of quotients for associative rings
given by Utumi in~\cite{utu}, the second author introduced
in~\cite{msm} the notion of algebra of quotients of a Lie algebra.
We now recall the main definitions and some results.

Let $L\subseteq Q$ be an extension of Lie algebras. For any $q$ in
$Q$, denote by $_L(q)$ the linear span in $Q$ of $q$ and the
elements of the form $\ad {x_1} \cdots \ad {x_n} q,$ where $n\in
\mathbb{N}$ and $x_1,\dots,  x_n \in L$. In particular, if $q\in
L$, then $_L(q)$ is just the ideal of $L$ generated by $q$.

We say that $Q$ is \emph{ideally absorbed into} $L$ if for every
nonzero element $q$ in $Q$ there exists an ideal $I$ of $L$ with
$\ann L(I)=0$ such that $0\neq [I,\ q]\subseteq L$.

\begin{definitions}{\rm (\cite[Definition 2.1 and Proposition
2.15]{msm}.) Let $L\subseteq Q$ be an extension of Lie algebras.
We say that $Q$ is an \emph{algebra of quotients of }$L$ (or also
that $L$ is a \emph{subalgebra of quotients of} $Q$) if the
following equivalent conditions are satisfied:
\begin{enumerate}[(i)]\itemsep=2mm
\item Given $p$ and $q$ in $Q$ with $p\neq 0$, there exists $x$ in $L$ such
that $[x, p]\neq 0$ and $[x,\ _L(q)]\subseteq L$.
\item $Q$ is ideally absorbed into $L$.
\end{enumerate}

If given a non-zero element $q$ in $Q$ there exists $x$ in $L$
such that $0\neq [x, q]\in L$, then $Q$ is said to be a
\emph{\wqa} of $L$.}
\end{definitions}
A Lie algebra $L$  has an \qa if and only if it has no nonzero
total zero divisors, or, equivalently, $\ann{}(L)=0$
(see~\cite[Remark 2.3]{msm}).
\begin{rem}
\label{29/03/04,1} If $Q$ is a \wqa of $L$, then $Z(Q)=Z(L)=0$.
Indeed, given a non-zero element $q$ in $Q$, there exists $x$ in
$L$ such that $0\neq [x, q]\in L$. Then $q\notin Z(Q)$. If $q\in
L$, then the same argument shows that $q\notin Z(L)$ and so $Z(Q)$
and $Z(L)$ have no non-zero elements.
\end{rem}

\begin{lem}
\label{29/03/04,2} Let $L\subseteq Q$ be an extension of Lie
algebras such that $Q$ is a \wqa of $L$, and let $I$ be an ideal
of $L$. If $\ann{L}(I)=0$ then $\ran{A(Q)}(A_Q(I))=0$.
\end{lem}
\begin{proof} Assume that $\ann{L}(I)=0$. We first note that $\ann Q
(I)=0$ since $Q$ is a \wqa of $L$ (and applying~\cite[Lemma 2.11
]{msm}). Now let $\mu\in \ran{A(Q)}(A_Q(I))$. Then $\ad y\mu=0$
for every $y$ in $I$. If $q\in Q$, we then have that $0=\ad
y\mu(q)=[y,\mu(q)]$. This says that $\mu (q)\in \ann{Q}(I)=0$, and
so $\mu=0$.

\end{proof}
\begin{lem}
\label{ayudaF1}Let $L\subseteq Q$ be an extension of Lie algebras,
and let $x_1,\dots, x_n, y\in L$. Then we have, in $A(Q)$:
\[
\ad{x_1}\cdots \ad{x_n}\ad y= \ad y \ad{x_1}\cdots
\ad{x_n}+{\sum_{i=1}^n} \ad{x_1}\cdots \ad{[x_i, y]}\cdots \
ad{x_n}\,.
\]
In particular, if $I$ is an ideal of $L$ and $x_1,\dots, x_n \in
I$, then
\[
\ad{x_1}\cdots \ad{x_n}\ad y = \ad y \ad{x_1}\cdots
\ad{x_n}+\alpha\,,
\]
where $\alpha\in \mathrm{span}\{\ad{z_1}\cdots \ad{z_{n}} \vert
z_i\in I\}$.
\end{lem}

\begin{proof}
The second part of the conclusion follows immediately from the
first. For this one, we use induction on $n$, the case $n=1$ being
obvious. If $n\geq 2$, then
\[
\ad{x_1}\cdots \ad{x_n}\ad y = \ad{x_1}\cdots
\ad{x_{n-1}}\ad{[x_n,y]}+ \ad{x_1}\cdots \ad{x_{n-1}}\ad y \ad
{x_n}\] which, by the induction step, is equal to
\[
\ad{x_1}\cdots \ad{x_{n-1}}\ad{[x_n,y]}+ {\sum_{i=1}^{n-1}}
\ad{x_1}\cdots \ad{[x_i, y]}\cdots \ad{x_{n-1}}\ad{x_n}+ (\ad y
\ad{x_1}\cdots \ad{x_{n-1}})\ad {x_n}\,,
\]
as wanted.
\end{proof}
Let $L\subseteq Q$ be an extension of Lie algebras. Denote by
$A_0$ be associative subalgebra of $A(Q)$ whose elements are those
$\mu$ in $A(Q)$ such that $\mu(L)\subseteq L$.  We obviously have
the containments:
\[
A_Q(L)\subseteq A_0 \subseteq A(Q)\,.
\]

In order to ease the notation in the next few results we will use
$\widetilde I$ to denote, for any ideal $I$ of $L$, the two-sided
ideal of $A_Q(L)$ generated by the elements of the form
$\ad{x}\colon Q\to Q$ for $x$ in $I$, i.e.
$\widetilde{I}=\mathrm{id}_{A_Q(L)}(A_Q(I))$.
\begin{lem}
\label{F1} Let $L\subseteq Q$ be an extension of Lie algebras. Let
$I$ be an ideal of $L$ and $q_1,\ldots q_n \in Q$ such that $[q_i,
I]\subseteq L$ for every $i=1, \ldots, n$. Then, for $\mu=\ad{q_1}
\cdots \ad{q_n}$ in $A(Q)$, we have that
$\mu\cdot\big(\widetilde{I}\big)^n\subseteq A_0$ and
$\big(\widetilde{I}\big)^n\cdot\mu\subseteq A_0$ (where
$\big(\widetilde{I}\big)^n$ denotes the $n$-th power of
$\widetilde{I}$ in the associative algebra $A_Q(L)$).
\end{lem}

\begin{proof}
Arguing as in Lemma~\ref{26/11/04,1}, it is enough to consider an
element $y$ in $\widetilde{I^n}$ of the form $y=\ad {x_1}\cdots
\ad{x_{n-1}}\ad {x_n}$, where $x_i\in I$, and prove that both $\mu
y$ and $y\mu$ belong to $A_0$. We will use induction on $n$. For
$n=1$ we have $\ad x\ad q= \ad{[x, q]}+\ad q\ad x$ and since $[x,
q]\in L$ we see that $\ad{[x,q]}\in A_0$. On the other hand, $\ad
q\ad x(L)\subseteq \ad q(I)\subseteq L$, and so $\ad q\ad x\in
A_0$.

Assume the result true for $n-1$. Now, by Lemma~\ref{ayudaF1} we
have
\begin{equation}
\label{igualdad}
\begin{split}
& \ad{x_1}\cdots\ad{x_n}\ad{q_1}\cdots\ad{q_n}=
(\ad{q_1}\ad{x_1})(\ad{x_2}\cdots\ad{x_n}\ad{q_2}\cdots\ad{q_n})+\\
& {\sum_{i=1}^n}\ad{x_1}\cdots\ad{[x_i,
q_1]}\cdots\ad{x_n}\ad{q_2}\cdots\ad{q_n}\,.
\end{split}
\tag{$\dagger$}
\end{equation}



The first summand in the last equality belongs to $A_0$ because,
as proved before, $\ad{q_1}\ad{x_1}\in A_0$ and
$\ad{x_2}\cdots\ad{x_n}\ad{q_2}\cdots\ad{q_n}\in A_0$ by the
induction hypothesis. On the other hand, for each of the terms
$\ad{x_1}\cdots\ad{[x_i,
q_1]}\cdots\ad{x_n}\ad{q_2}\cdots\ad{q_n}$ we have that $x_i\in I$
and $[x_i, q_1]\in L$. Using Lemma~\ref{ayudaF1} we may write this
as:
\[
\ad{[x_i,q_1]}\ad{x_1}\cdots\ad{x_{i-1}}\ad{x_{i+1}}\cdots\ad{x_n}\ad{q_2}\cdots\ad{q_n}+
\alpha\cdot \ad{x_{i+1}}\cdots\ad{x_n}\ad{q_2}\cdots\ad{q_n}\,,
\]
where $\alpha \in \mathrm{span}\{\ad{z_1}\cdots\ad{z_{i-1}} \vert
z_j\in I\}$. The induction hypothesis applies again to show that
this belongs to $A_0$. Hence $y\mu\in A_0$.

If we continue to develop in the expression (\ref{igualdad}), we
get that, for some $\alpha_0$ in $A_0$,
\[
\ad{q_1}\ad{q_2}\ad{x_1}\cdots\ad{x_n}\ad{q_3}\cdots\ad{q_n}+
{\sum_{i=1}^n}\ad{q_1}\ad{x_1}\cdots\ad{[x_i,
q_2]}\cdots\ad{x_n}\ad{q_3}\cdots\ad{q_n}+ \alpha_0\,.
\]
Using Lemma~\ref{ayudaF1} we can write each term of the form
\[
\ad{q_1}\ad{x_1}\cdots\ad{[x_i,
q_2]}\cdots\ad{x_n}\ad{q_3}\cdots\ad{q_n}
\]
as:
\[
\ad{q_1}\ad{[x_i,
q_2]}\ad{x_1}\cdots\ad{x_{i-1}}\ad{x_{i+1}}\cdots\ad{x_n}\ad{q_3}\cdots\ad{q_n}+
\ad{q_1}\cdot\alpha\cdot\ad{x_{i+1}}\cdots\ad{x_n}\ad{q_3}\cdots\ad{q_n}\,,
\]
where $\alpha\in \mathrm{span}\{\ad{z_1}\cdots\ad{z_{i-1}}\vert
z_j\in I\}$. Since $[x_i, q_2]\in L$ and $x_i\in I$, we have that
$\ad{q_1}\ad{[x_i, q_2]}\ad{x_1}=\ad{q_1}\ad{[[x_i, q_2], x_1]}+
\ad{q_1}\ad{x_1}\ad{[x_i, q_2]}$. From this we see that the first
summand above belongs to $A_0$. For the second summand, assuming
that $\alpha=\ad{z_1}\cdots\ad{z_{i-1}}$, with $z_j$ in $I$, we
have
$(\ad{q_1}\ad{z_1})\ad{z_2}\cdots\ad{z_{i-1}}\ad{x_{i+1}}\cdots\ad{x_n}\ad{q_3}\cdots\ad{q_n}$,
which is also an element of $A_0$. Continuing in this way, we find
that
\[
\ad{x_1}\dots\ad{x_n}\ad{q_1}\dots\ad{q_n}-\ad{q_1}\dots\ad{q_n}\ad{x_1}\dots\ad{x_n}\in
A_0
\]
and by what we have just proved, we see that
$\ad{q_1}\dots\ad{q_n}\ad{x_1}\dots\ad{x_n}\in A_0$, as was to be
shown.
\end{proof}

\begin{cor}
\label{cordeF1} Let $L\subseteq Q$ be an extension of Lie
algebras. Let $\mu=\ad{q_1}\cdots \ad{q_n}$ be in $A(Q)$ and
suppose  that there exists an  ideal $I$ of $L$ that satisfies
$[q_i, I]\subseteq L$ for every $i=1, \dots, n$. Then
$\mu\widetilde{I^n}\subseteq A_0$ and $\widetilde{I^n}\mu\subseteq
A_0$ (where $I^n$ denotes the $n$-th power of I in the Lie algebra
$L$).
\end{cor}

\begin{lem}
\label{facil}

Let $L$ be a semiprime Lie algebra. If $I$ is an ideal of $L$ with
$\ann L (I)=0$, then $\ann L(I^s)=0$ for any $s\geq 1$. Any
intersection of ideals with zero annihilator will also have zero
annihilator.
\end{lem}

\begin{proof}
For semiprime Lie algebras, ideals with zero annhilator are
exactly the essential ideals, by~\cite[Lemma 1.2 (ii)]{msm}. Hence
the conclusion follows at once.
\end{proof}

\begin{thm}
\label{elteorema} Let $L\subseteq Q$ be an extension of Lie
algebras with $L$ semiprime. Then the following conditions are
equivalent:
\begin{enumerate}[{\rm (i)}]\itemsep=2mm
\item $Q$ is an \qa of $L$,
\item $Z(Q)=0$ and, if $\mu \in A(Q)\setminus \{0\}$, there is
an ideal $I$ of $L$ with $\ann L(I)=0$ such that
$\mu\widetilde{I}\subseteq A_0$, $0\neq \widetilde{I}\mu \subseteq
A_0$. If $\mu=\ad{q}$, then we also have $\mu\widetilde{I}(L)\neq
0$.
\end{enumerate}
\end{thm}
\begin{proof}

(ii)$\Rightarrow$ (i)\,: Let $q\in Q\smallsetminus\{0\}$. Then
$\mu=\ad{q}\neq 0$ since $Z(Q)= 0$. Let $\widetilde{I}$ be as in
(ii), so it satisfies $\mu\widetilde{I}(L)\neq 0$ and
$\mu\widetilde{I}\subseteq A_0$. Set
\[
I_0:=\mathrm{span}\{\alpha(x) \ \vert \ x\in L\text{ and } \alpha
\in \widetilde{I}\}\,.
\]
Then $I_0$ is an ideal of $L$ such that $\ann L(I_0)=0$. Indeed,
if $x$, $y\in L$ and $\alpha \in\widetilde I$, we have $[y,
\alpha(x)]= (\ad{y}\alpha)(x)$ and $\ad{y}\alpha\in
\widetilde{I}$. If now $[x, I_0 ]=0$ for some $x$ in $L$, then
$\ad{x}\widetilde{I}(L)=0$. In particular, for $y$ and $z$ in $I$
we have that $\ad{x}\ad{y}(z)=0$, so $x\in \ann L([I, I])$, which
is zero by Lemma~\ref{facil}.

Finally, $0\neq [q, I_0]=\ad{q}\widetilde{I}(L)\subseteq L$.

(i)$\Rightarrow$(ii)\,: Since $Q$ is an algebra of quotients of
$L$, it is also a weak algebra of quotients of $L$, hence
Remark~\ref{29/03/04,1} applies in order to obtain that $Z(Q)=0$.

Next, let $\mu=\sum_{i\geq 1}\ad{q_{i,1}}\cdots \ad{q_{i, r_i}}\in
A(Q)\setminus\{0\}$. Using that $Z(Q)=0$ we may of course assume
that all $q_{i,j}$ are non-zero elements in $Q$. Set
$s=\sum_{i\geq 1}r_i$. As $Q$ is an \qa of $L$, there exists, for
every $i$ and $j$, an ideal $J_{i, j}$ of $L$ such that
$\ann{L}(J_{i, j})=0$ and $0\neq [J_{i, j}, q_{i, j}]\subseteq L$.
By Lemma~\ref{facil}, the ideal $J=\bigcap_{i, j}J_{i,j}$ and
hence also $I=J^s$ will have zero annihilator in $L$. Then $[q_{i,
j}, I]\subseteq [q_{i, j}, J_{i, j}]\subseteq L$. By
Corollary~\ref{cordeF1} and taking into account that $J^s\subseteq
J^{r_i}$ for every $i$, we have that $\mu \widetilde{I}\subseteq
A_0$ and $ \widetilde{I}\mu\subseteq A_0$.

If $\widetilde{I}\mu=0$, then $\mu\in\ran{A(Q)}(A_Q(I))$ which is
zero by Lemma~\ref{29/03/04,2}.

Finally, suppose that $\mu=\ad q$, where $q\in
Q\smallsetminus\{0\}$. By \cite[Lemma 2.13]{msm}, $Q$ is an \qa of
$I^2$. This implies that there exist elements $y$ and $z$ in $I$
such that $0\neq [q, [y, z]]$. But this means that $\ad q\ad y
z\in \ad q{\widetilde{I}}(L)$.
\end{proof}

Recall that an associative algebra $S$ is said to be a \emph{left
quotient algebra} of a subalgebra $A$ if whenever $p$ and $q\in S$
with $p\neq 0$, there exists $x$ in $A$ such that $xp\neq 0$ and
$xq\in A$. An associative algebra $A$ has a left quotient algebra
if and only if it has no total right zero divisors different from
zero. (Here, an element $x$ in $A$ is said to be a \emph{total
right zero divisor} if $Ax=0$.)
\begin{cor}
\label{leftquotient} Let $L\subseteq Q$ be an extension of Lie
algebras with $L$ semiprime. Suppose that $Q$ is an \qa of $L$.
Then $A(Q)$ is a left quotient algebra of $A_0$.
\end{cor}

\begin{proof} We will use the characterization of left quotient
(associative) algebras given in~\cite[Lemma 2.14]{msm}. Let
$\mu\in A(Q)\smallsetminus\{0\}$ and let $I$ be an ideal of $L$
satisfying condition (ii) in Theorem~\ref{elteorema}. Set
$J=A_0\widetilde{I}+\widetilde{I}$, a left ideal  of $A_0$ that
satisfies $0\neq J\mu \subseteq A_0$ (because $0\neq
\widetilde{I}\mu \subseteq A_0$).

Since also $\ann L(I)=0$ we obtain from Lemma~\ref{29/03/04,2}
that $\ran{A(Q)}(A_Q(I))=0$. This, together with the fact that
$\widetilde{I}=\idl{A_Q(L)}(A_Q(I))=\idr {A_Q(L)}(A_Q(I))=\id
{A_Q(L)}(A_Q(I))$ (see Lemma~\ref{26/11/04,1}), implies that
$\ran{A(Q)}(\widetilde{I})=0$. Since $\ran{A_0}(J)\subseteq
\ran{A(Q)}(\widetilde I)$ we get that also $\ran{A_0}(J)=0$. This
concludes the proof.
\end{proof}

\begin{rem}
\label{obs} Note that in the situation of Theorem~\ref{elteorema}
the ideal $\widetilde{I}$ has zero right annihilator in $A(Q)$ and
therefore, as is established in the proof of the corollary above,
$A_0\widetilde{I}+\widetilde{I}$ is a left ideal of $A_0$ with
zero right annhilator.
\end{rem}

In a somewhat different direction, we analyse other instances
where the notion of an algebra of quotients in the associative
case is closely related to the one in the Lie case. If $A$ is any
semiprime (associative) algebra, we denote by
$Q_{\mathrm{max}}^l(A)$ the maximal left quotient algebra of $A$
and by $Q_s(A)$ the Martindale symmetric algebra of quotients of
$A$. Recall that $Q_s(A)$ can be characterized as those elements
$q$ in $Q_{\mathrm{max}}^l(A)$ for which there is an essential
ideal $I$ of $A$ satisfying $Iq+qI\subseteq A$ (see,
e.g.~\cite[Section 2.2]{bmm}).

We remind the reader that for an associative algebra $A$ we denote
by $A^{(-)}$ the Lie algebra which equals $A$ as a module and has
Lie bracket given by $[a,b]=ab-ba$.

It was proved by the second author in~\cite[Proposition 2.16]{msm}
that, if $A\subseteq Q\subseteq Q_s(A)$ are associative algebras
with $A$ semiprime and $Z(Q)=0$, then $Q^{(-)}$ is an algebra of
quotients of $A^{(-)}$. With considerably more effort we shall
establish an improvement of this result that provides with new
examples of algebras of quotients of Lie algebras.

In order to prove Theorem~\ref{liecenter} below, we need a lemma
which is an adaptation of~\cite[Lemma 2]{hers} to our setting. We
include the necessary changes in the statement and proof. Recall
that an algebra $A$ is \emph{$2$-torsion free} provided that
$2x=0$ implies $x=0$.
\begin{lem}
\label{herstein} Let $A$ be a semiprime $2$-torsion free algebra
and let $Q$ be a subalgebra of $Q_s(A)$ that contains $A$. Let
$q\in Q$ and assume that there is an essential ideal $I$ of $A$
such that $qI+Iq\subseteq A$ and $[q,[I,I]]=0$. Then $[q,I]=0$.
\end{lem}

\begin{proof}
We first note that $Q$ is also $2$-torsion free. For, if $q$ is a
non-zero element in $Q$ and $2q=0$, then, since $Q$ is an algebra
of quotients of $A$, there exists $a$ in $A$ such that $aq$ is a
non-zero element of $A$. But $2(aq)=0$, which contradicts our
assumption that $A$ is $2$-torsion free.

Next, define $d\colon A\to Q$ by $d(x)=[x,q]$. As in the proof
of~\cite[Lemma 2]{hers}, one verifies using $2$-torsion freeness
of $Q$ that
\[
d(x)(uv-vu)=0\,,\,\text{whenever }x\in A\,,\,\text{and }u,v\in
[I,I]\,.
\]
It also follows as in~\cite[Lemma 2]{hers} that the set
\[
J=\{r\in A\mid d(x)r=0\text{ for all }x\text{ in }A\}
\]
is an ideal of $A$ and we have $[[I,I],[I,I]]\subseteq J$.

We now proceed to prove that $A/J$ is semiprime and $2$-torsion
free. Denote by $\overline{x}$ the class modulo $J$ of an element
$x$ in $A$. If $2\overline{x}=0$, then $2x\in J$ and hence
$d(y)(2x)=0$ for all $y$ in $A$. Since $Q$ is $2$-torsion free
this implies that $x\in J$, that is, $\overline{x}=0$.

Now, let $K$ be an ideal of $A$ (containing $J$) for which
$\overline{K}^2=0$, that is, $K^2\subseteq J$. We then have that
$d(y)K^2=0$ for any $y$ in $A$.

By assumption, $Iq+qI\subseteq A$, and so $Id(y)\subseteq A$ and
also $d(y)KI$ is a right ideal of $A$ for any $y$ in $A$. Then
$(d(y)KI)^2=d(y)KId(y)KI\subseteq d(y)K^2=0$ and it follows by
semiprimeness of $A$ that $d(y)KI=0$. Since $I$ is an essential
ideal in $A$ and $Q$ is an algebra of quotients of $A$, $I$ has
zero left annihilator in $Q$ and therefore $d(y)K=0$. By the
definition of $J$, this shows that $K\subseteq J$, i.e.
$\overline{K}=0$, and thus $A/J$ is semiprime.

The ideal $\overline{I}:=(I+J)/J$ of $A/J$ clearly satisfies
$[[\overline{I},\overline{I}],[\overline{I},\overline{I}]]=0$.
Using~\cite[Lemma 1]{hers} twice we conclude that
$\overline{I}\subseteq Z(A/J)$, that is,
$[\overline{I},\overline{A}]=0$. Thus $[I,A]\subseteq J$ and
therefore $d(y)[I,A]=0$ for all $y$ in $A$.

Let $K_1=\{x\in A\mid x[I,A]=0\}$, which is an ideal of $A$. By
what we have just proved, $d(u)\in K_1$ whenever $u\in A$.

We now claim that $K_1I^2d(u)=0$ for any $u$ in $I$. To see this,
take $x$ in $K_1$, $v$, $w$ in $I$ and compute that (using
$Iq+qI\subseteq A$):
\begin{align*}
xvwd(u) &=xvw(uq-qu) =x(vwuq-uqvw+uqvw-vwqu)\\
&=x[vw,uq]+x(uqvw-vwuq) =0+x(uqvw-qvwu+qvwu-vwqu)\\
&=x[u,qvw]+x(qvw-vwq)u=0+x(qvw-wqv+wqv-vwq)u\\
&=x([qv,w]+[wq,v])u=0\,,
\end{align*}
from which the claim follows.

Since $I^2$ is also essential in $A$ we obtain that $d(u)K_1=0$.
Altogether this implies that $d(u)\in K_1\cap
\mathrm{l.ann}_A(K_1)$, which is zero as it is a nilpotent left
ideal in $A$. Thus $d(u)=0$ for all $u$ in $I$, that is,
$[q,I]=0$, as desired.
\end{proof}
\begin{thm}
\label{liecenter} Let $A$ be a semiprime $2$-torsion free
associative algebra and let $Q$ be a subalgebra of $Q_s(A)$ that
contains $A$. Then $A^{(-)}/Z(A)$ and $Q^{(-)}/Z(Q)$ are semiprime
Lie algebras and $Q^{(-)}/Z(Q)$ is a (Lie) algebra of quotients of
$A^{(-)}/Z(A)$.
\end{thm}

\begin{proof}
Let $\overline{K}$ be a Lie ideal of $A^{(-)}/Z(A)$ such that
$[\overline{K},\overline{K}]=0$. Then $\overline{K}$ is the image
of a Lie ideal $K$ of $A^{(-)}$ via the natural map $A^{(-)}\to
A^{(-)}/Z(A)$. The condition on $\overline{K}$ translates upstairs
into $[K,K]\subseteq Z(A)$. By~\cite[Lemma 1]{hers} we have
$K\subseteq Z(A)$, that is, $\overline{K}=0$.

Since $Q$ is also semiprime (e.g.~\cite[Lemma 2.1.9 (i)]{bmm}) and
$2$-torsion free (by the first part of the proof of
Lemma~\ref{herstein}) the same argument applies to show that
$Q^{(-)}/Z(Q)$ is semiprime.

By~\cite[Lemma 1.3 (i)]{msm}, $Z(A)=Z(Q)\cap A$. Therefore the
natural map
\[
A^{(-)}/Z(A)\to Q^{(-)}/Z(Q)
\]
is an inclusion and we shall identify $A^{(-)}/Z(A)$ with its
image into $Q^{(-)}/Z(Q)$. We now prove that $Q^{(-)}/Z(Q)$ is
ideally absorbed into $A^{(-)}/Z(A)$.

For an element $x$ in $Q$, denote by $\overline{x}$ the class of
$x$ in $Q^{(-)}/Z(Q)$. Let $\overline{q}$ be a non-zero element in
$Q^{(-)}/Z(Q)$. Since by assumption $Q\subseteq Q_s(A)$ there
exists an essential ideal $I$ of $A$ such that $Iq+qI\subseteq A$.

We claim that the Lie ideal $\overline{I}=(I^{(-)}+Z(A))/Z(A)$ has
zero annihilator in $A^{(-)}/Z(A)$. For $x$ in $A$, if
$\overline{x}$ satisfies $[\overline{x},\overline{I}]=0$, then
$[x,I^{(-)}]=[x,I]\subseteq Z(A)$. Therefore $[[x,I],I]=0$.
Applying the Jacobi identity this yields $[x,[I,I]]=0$.
Therefore~\cite[Lemma 2]{hers} allows to conclude that $[x,I]=0$.
Note that $A$ is an algebra of quotients of $I$ since it is an
essential ideal of $A$. Hence, using~\cite[Lemma 1.3 (iv)]{msm} we
obtain $[x,A]=0$. Thus $\overline{x}=0$ in $A^{(-)}/Z(A)$.

Next, $[\overline{I},\overline{q}]=\overline{[I^{(-)},q]}\subseteq
A^{(-)}/Z(A)$ because $Iq+qI\subseteq A$.

Finally, we need to see that $[\overline{I}, \overline{q}]\neq 0$.
To reach a contradiction, suppose that $[I,q]\subseteq Z(A)$. Then
$[[q,I],I]=0$ and using the Jacobi identity we have $[q,[I,I]]=0$.
By Lemma~\ref{herstein} this implies $[q,I]=0$, and a second use
of~\cite[Lemma 1.3 (iv)]{msm} yields $[q,Q]=0$, that is,
$\overline{q}=0$, which contradicts the choice of $q$.
\end{proof}

We will now draw a consequence of Theorem~\ref{liecenter}. First
we need a lemma.

\begin{lem}
\label{quotientcomm} Let $L\subseteq Q$ be an extension of Lie
algebras with $L$ semiprime and such that $Q$ is an algebra of
quotients of $L$. Then $[Q,Q]$ is an algebra of quotients of
$[L,L]$.
\end{lem}
\begin{proof}
Let $\sum\limits_{i=1}^n [x_i,y_i]$ be a non-zero element in
$[Q,Q]$. Since $Q$ is an algebra of quotients of $L$ and $L$ is
semiprime, we may choose an ideal $I$ of $L$ with $\ann L(I)=0$
such that $[x_i,I]$, $[y_i,I]$, $[[x_i,y_i],I]\subseteq L$ for all
$i$.

We know that $[I,I]$ also has zero annihilator in $L$, and from
the inclusion
\[
\ann {[L,L]}([I,I])\subseteq \ann L ([I,I])
\]
we obtain $\ann {[L,L]}([I,I])=0$. Hence, $\ann Q([I,I])=0$
using~\cite[Lemma 2.11]{msm}. In particular, $[\sum\limits_{i=1}^n
[x_i,y_i],[I,I]]\neq 0$.

Moreover, for every $z$, $t$ in $I$ and each $i$, we have that
\[
[[x_i,y_i],[z,t]]=[[[x_i,y_i],z],t]+[z,[[x_i,y_i],t]]\in [L,L]\,,
\]
using the Jacobi identity.
\end{proof}

\begin{cor}
\label{liequocenter} Let $A$ be a semiprime $2$-torsion free
associative algebra and let $Q$ be a subalgebra of $Q_s(A)$ that
contains $A$. Then $[Q^{(-)},Q^{(-)}]/Z([Q^{(-)},Q^{(-)}])$ is a
quotient algebra of $[A^{(-)},A^{(-)}]/Z([A^{(-)},A^{(-)}])$.
\end{cor}
\begin{proof}
First note that $Z(A^{(-)})=Z(A)=Z([A^{(-)},A^{(-)}])$
by~\cite[Lemma 2]{hers}, and the same conclusion holds for $Q$.
From this it follows that
$[A^{(-)}/Z(A),A^{(-)}/Z(A)]=[A^{(-)},A^{(-)}]/Z([A^{(-)},A^{(-)}])$,
and analogously for $Q$. The result is then obtained applying
Theorem~\ref{liecenter} and Lemma~\ref{quotientcomm}.
\end{proof}

\section{Multiplicatively semiprime Lie algebras and dense extensions}

Given an extension $L\subseteq Q$ of Lie algebras, we have
considered in the previous section the associative algebra $A(Q)$
and the subalgebra $A_Q(L)$ generated by the elements $\ad x$ for
$x$ in $L$. It is natural to study the relationship between this
algebra and the associative algebra $A(L)$.

Given an element $\mu$ in $A(L)$, we can of course think of this
element in $A(Q)$ because of the inclusion $L\subseteq Q$. In
order to distinguish this change of domains, we will use the
notation $\mu^L$ and $\mu^Q$. Thus, for example, given $x$ in $L$
we have $\ad{x}^L$ and $\ad{x}^Q$ which differ in the use of the
bracket in $L$ and in $Q$ respectively. With these considerations
in mind it is desirable to have a well-defined map $\varphi\colon
A(L)\to A(Q)$ given by $\mu^L\mapsto \mu^Q$. Whilst it is not
guaranteed this can be done, if such a map exists then it is an
injective algebra homomorphism and $\varphi(A(L))=A_Q(L)$.

The condition expressed above is just a rephrasing of the density
condition introduced by Cabrera in~\cite{Ca}. Specifically, for
any algebra $L$ (not necessarily associative) over our commutative
ring of scalars $\Phi$, let $M(L)$ be the subalgebra of
$\mathrm{End}_{\Phi}(L)$ generated by the identity map together
with the operators given by right and left multiplication by
elements of $L$. In the case of a Lie algebra $L$, then $M(L)$ is
nothing but the unitization of $A(L)$.

Following~\cite{Ca}, given an extension of (not necessarily
associative) algebras $L\subseteq Q$, the \emph{annihilator of}
$L$ \emph{in} $M(Q)$ is defined by:
\[
L^{ann}:=\{ \mu \in M(Q)\vert \ \mu(x)=0 \hbox{ for every } x \in
\ L\}\,.
\]
If $L^{ann}=0$, then $L$ is said to be a \emph{dense subalgebra}
of $Q$, and we will say that $L\subseteq Q$ is a \emph{dense
extension of algebras}.
\begin{lem}

\label{densext} Let $L\subseteq Q$ be a dense extension of Lie
algebras. If $Z(Q)=0$ then $Z(L)=0$.
\end{lem}

\begin{proof}
Suppose that $x\in L$ satisfies $[x,L]=0$. This means that $\ad
x(L)=0$, and hence $\ad x=0$ as an element of $A(Q)$. Since
$Z(Q)=0$, it follows that $x=0$.
\end{proof}

\begin{lem}
\label{20/02/05,1} Let $L\subseteq Q$ be an extension of Lie
algebras and suppose that $Z(Q)=0$. Then the following conditions
are equivalent:
\begin{enumerate}[{\rm (i)}]\itemsep=2mm
\item $L$ is a dense subalgebra of $Q$.
\item If $\mu(L)=0$ for some $\mu$ in $A(Q)$, then $\mu=0$.
\end{enumerate}

\end{lem}
\begin{proof} Clearly, (i) implies (ii). Conversely, suppose that
$\mu\in M(Q)$ satisfies $\mu(L)=0$. If $\mu(p)\neq 0$ for some $p$
in $Q$, then use $Z(Q)=0$ to find a (non-zero) element $q$ in $Q$
satisfying $\ad q \mu(p)\neq 0$. But then $\ad q \mu(L)= 0$ and
since $A(Q)$ is a two-sided ideal of $M(Q)$, we have that $\ad
q\mu\in A(Q)$. Hence condition (ii) yields $\ad q \mu = 0$, a
contradiction.
\end{proof}

We will now present some examples of dense extensions.

\begin{lem}
\label{densecomm} Let $L\subseteq Q$ be a dense extension of Lie
algebras. Then $[L,L]\subseteq [Q,Q]$ is also dense.
\end{lem}

\begin{proof}
Let $\mu\in M([Q,Q])$ and suppose that $\mu ([L,L])=0$. Then, for
any $x$ in $L$, we have $\mu\ad x(L)=0$. Since $\mu\ad x\in M(Q)$
and $L\subseteq Q$ is a dense extension we obtain $\mu\ad x(Q)=0$.
This implies that $\mu\ad q (x)=-\mu\ad x(q)=0$ for every $q$ in
$Q$ and every $x$ in $L$. A second use of density implies that
$\mu\ad q(Q)=0$, that is, $\mu ([Q,Q])=0$.
\end{proof}

\begin{lem}
\label{denseassoc} Let $A$ be a semiprime associative algebra.
Then, for any subalgebra $Q$ of $Q_{\mathrm{max}}^l(A)$ that
contains $A$, the extension $A\subseteq Q$ is dense.
\end{lem}

\begin{proof}
Let $\mu\in M(Q)$ such that $\mu(A)=0$. We may write
$\mu(x)=\sum\limits_{i=1}^n p_ixq_i$ with $p_i$, $q_i$ in
$Q_{\mathrm{max}}^l(A)$. Then $\mu(x)$ is a Generalized Polynomial
Identity on $A$ (in the terminology of~\cite[p. 212]{bmm}), and
so~\cite[Proposition 6.3.13]{bmm} applies to conclude
$\mu(Q_{\mathrm{max}}^l(A))=0$. Therefore $\mu(Q)=0$.
\end{proof}

In the following lemma we consider an extension $A\subseteq Q$ of
associative algebras where $Q$ is a left quotient algebra of $A$.
We then have that $Z(A)=Z(Q)\cap A$ (\cite[Lemma 1.3 (i)]{msm}) so
that we may consider, as in Theorem~\ref{liecenter}, that
$A^{(-)}/Z(A)$ is a subalgebra of $Q^{(-)}/Z(Q)$.
\begin{prop}
\label{denselie} Let $A$ be an associative algebra and $Q$ a
subalgebra of $Q_{\mathrm{max}}^l(A)$ that contains $A$. Then
$A^{(-)}/Z(A)\subseteq Q^{(-)}/Z(Q)$ and
\[
[A^{(-)},A^{(-)}]/Z([A^{(-)},A^{(-)}])\subseteq
[Q^{(-)},Q^{(-)}]/Z([Q^{(-)},Q^{(-)}])
\]
are dense extensions of Lie algebras.
\end{prop}

\begin{proof}
Denote by $\overline{q}$ the classes of elements in
$Q^{(-)}/Z(Q)$. Let $\mu\in M(Q^{(-)}/Z(Q))$ such that
$\mu(A^{(-)}/Z(A))=0$. Then there is $\mu'$ in $M(Q)$ such that
$\overline{\mu'(q)}=\mu(\overline{q})$ for every $q$ in $Q$.

The condition that $\mu$ vanishes on $A^{(-)}/Z(A)$ means that
$\mu'(A)\subseteq Z(Q)$. Suppose that $\mu'(Q)$ is not contained
in $Z(Q)$. Then there are non-zero elements $p$ and $q$ in $Q$ for
which $(\lambda_q\circ\mu'-\rho_q\circ\mu')(p)=[q,\mu'(p)]\neq 0$,
where $\lambda_q$ and $\rho_q$ stand for the left and right
(associative) multiplication by $q$, respectively.

But then $\lambda_q\circ\mu'-\rho_q\circ\mu'$ is a non-zero
element in $M(Q)$ and since the extension $A\subseteq Q$ is dense
by Lemma~\ref{denseassoc} we get that
$[q,\mu'(A)]=(\lambda_q\circ\mu'-\rho_q\circ\mu')(A)\neq 0$. This
contradicts the fact that $\mu'(A)\subseteq Z(Q)$. This shows
$A^{(-)}/Z(A^{(-)})\subseteq Q^{(-)}/Z(Q^{(-)})$ is dense.

Finally, $[A^{(-)},A^{(-)}]/Z([A^{(-)},A^{(-)}])\subseteq
[Q^{(-)},Q^{(-)}]/Z([Q^{(-)},Q^{(-)}])$ is also dense by the first
part of the proof and Lemma~\ref{densecomm}.
\end{proof}

Following Cabrera and Mohammed~\cite{CabM}, we say that an algebra
$L$ is \emph{multiplicatively semiprime} (respectively,
\emph{prime}) whenever $L$ and its multiplication algebra $M(L)$
are semiprime (respectively, prime). Observe that in this
situation, and if $L$ is a Lie algebra, then $A(L)$ will also be a
semiprime (respectively, prime) algebra.

As we have seen, under mild hypotheses we have dense extensions of
semiprime Lie algebras 
$
A^{(-)}/Z(A)\subseteq Q^{(-)}/Z(Q)$ and
\[
[A^{(-)},A^{(-)}]/Z([A^{(-)},A^{(-)}])\subseteq
[Q^{(-)},Q^{(-)}]/Z([Q^{(-)},Q^{(-)}])\,,
\] 
so that $Q^{(-)}/Z(Q)$ is
a quotient algebra of $A^{(-)}/Z(A)$ and
$[Q^{(-)},Q^{(-)}]/Z([Q^{(-)},Q^{(-)}])$ is a quotient algebra of
$[A^{(-)},A^{(-)}]/Z([A^{(-)},A^{(-)}])$.

We remark that if $A$ is a semiprime algebra over a field of characteristic not $2$, then $A^{(-)}/Z(A)$ is multiplicatively semiprime (or even multiplicatively prime) in some important cases that are
covered in~\cite[Corollary 2.4]{CCLM}, but not in general
(see~\cite[Theorem 2.1]{CCLM}). This contrasts with the case of
$[A^{(-)},A^{(-)}]/Z([A^{(-)},A^{(-)}])$, which is always
multiplicatively semiprime if $A$ is semiprime (\cite[Corollary
2.4]{CCLM}).

If, furthermore, our algebra $A$ is endowed with an involution $^*$, denote by $K_A$ the set of all skew elements of $A$, that is, $K_A=\{x\in A\mid x^*=-x\}$. This is a Lie subalgebra of $A^{(-)}$, and it turns out that $K_A/Z(K_A)$ is multiplicatively semiprime or prime in various instances (see~\cite[Theorems 3.4, 3.6]{CCLM}). Again the situation is different for the Lie algebra $[K_A,K_A]/Z([K_A, K_A])$, which is always multiplicatively semiprime (\cite[Theorem 2.3]{CCLM}). The abovementioned examples involving commutators are of great interest
since they appear in Zelmanov's classification of simple
$M$-graded Lie algebras over a field of characteristic at least
$2d+1$, where $d$ is the width of $M$ (see~\cite{Zel}). It seems plausible (to the authors) that for algebras of the latter type -- $K_A/Z(K_A)$ in the case where $A$ has an involution -- results analogous
to Theorem~\ref{liecenter} and Proposition~\ref{denselie} are
available.

\begin{lem}
\label{20/02/05,2} Let $L\subseteq Q$ be a dense extension of Lie
algebras. If $Q$ is multiplicatively semiprime, then $A_0$ is
semiprime.
\end{lem}

\begin{proof}
Let $I$ be an ideal of $A_0$ with $I^2=0$. As in the proof of
Theorem~\ref{elteorema}, let $I_0=\mathrm{span}\{\alpha(x)| \alpha
\in I\,,\, x\in L\}$, which is clearly an ideal of $L$. For any
$\mu$ in $I$, we evidently have $\mu(I_0)=0$. It then follows
from~\cite[Proposition 3.1]{Ca} that $\mu(M(Q)(I_0))=0$ (where
$M(Q)(I_0)$ is the set of finite sums of elements of $M(Q)$
applied to elements of $I_0$). This implies that $\mu
M(Q)\mu(L)=0$, and thus $\mu M(Q)\mu=0$ since $L$ is dense in $Q$.
But $M(Q)$ is semiprime by hypothesis, so $\mu=0$ and since $\mu$
was an arbitrary element of $I$, we get that $I=0$, that is, $A_0$
is semiprime.
\end{proof}

\begin{rem}
\label{20/02/05,3} Given a dense extension of Lie algebras
$L\subseteq Q$, where $Q$ is multiplicatively semiprime, we have
that $L$ is also multiplicatively semiprime. This follows
from~\cite[Proposition 2.2]{Ca}.
\end{rem}

\begin{prop}
\label{20/02/05,5} Let $L\subseteq Q$ be a dense extension of Lie
algebras. Assume that $Q$ is a multiplicatively semiprime algebra
of quotients of $L$. Then $I\subseteq Q$ is a dense extension for
every essential ideal $I$ of $L$.
\end{prop}

\begin{proof}
We first observe that $Z(Q)=0$ because $Q$ is semiprime and hence
Lemma~\ref{20/02/05,1} applies. Thus, let $\mu$ be in $A(Q)$ such
that $\mu(I)=0$, and by way of contradiction assume that $\mu\neq
0$. By Corollary~\ref{leftquotient} (which can be applied by
virtue of Remark~\ref{20/02/05,3}) $A(Q)$ is a left quotient
algebra of $A_0$ and hence there exists $\lambda$ in $A_0$ such
that $0\neq \lambda \mu\in A_0$. Since the extension $L\subseteq
Q$ is dense, $\lambda\mu(L)\neq 0$, and since $\ann L(I)=0$, there
is a (non-zero) element $y$ in $I$ such that $\ad y
\lambda\mu(L)\neq 0$. Using now that $A(Q)$ has no total right
zero divisors, we get that $A(Q)\ad y\lambda \mu\neq 0$, and this,
coupled with the semiprimeness of $A(Q)$, implies that $A(Q)\ad y
\lambda\mu A(Q)\ad y \lambda\mu\neq 0$. A second application of
the fact that $L\subseteq Q$ is a dense extension yields $A(Q)\ad
y \lambda\mu A(Q)\ad y \lambda\mu(L)\neq 0$. However, $\mu(I)=0$
by assumption and thus~\cite[Proposition 3.1]{Ca} implies that
$\mu M(Q)(I)=0$. But this is a contradiction, because of the
containments
\[
\mu A(Q)\ad y \lambda\mu(L)\subseteq \mu A(Q)([I, L]) \subseteq
\mu A(Q)(I)\subseteq \mu M(Q)(I)=0\,.
\]
\end{proof}

\begin{cor}
{\rm (cf.~\cite[Lemma 3.1]{CCLM})} Let $L\subseteq Q$ be a dense
extension of Lie algebras with $Q$ multiplicatively prime. If $Q$
is a quotient algebra of $L$, then $I\subseteq Q$ is a dense
extension for any non-zero ideal $I$ of $L$.
\end{cor}

\begin{proof}
This follows directly from the previous proposition once we
realize that, in a prime Lie algebra, every non-zero ideal is
essential.
\end{proof}

\begin{cor}
\label{20/02/05,6} Let $L\subseteq Q$ be a dense extension of Lie
algebras. Assume that $Q$ is a multiplicatively semiprime algebra
of quotients of $L$. Then, for every essential ideal $I$ of $L$,
$\lan{A(Q)}(\widetilde{I})=0$.
\end{cor}

\begin{proof}
Let $\mu\in\lan{A(Q)} (\widetilde{I})$. Then, if $y\in I$, we have
$\mu\ad y(I)=0$. This implies that $\mu(I^2)=0$. By
Proposition~\ref{20/02/05,5} applied to the essential ideal $I^2$
of $L$, the extension $I^2\subseteq Q$ is dense, and so $\mu=0$.
\end{proof}

We now tie up our results with the associative case by relating
them to the maximal symmetric ring of quotients. This was first
introduced by Schelter in~\cite{sch} and systematically explored
by Lanning in~\cite{lann}. It has more recently come into
prominence following~\cite{edu}. In short, the maximal symmetric
ring of quotients $Q_{\sigma}(R)$ of a ring $R$ is the subring of
$Q_{\mathrm{max}}^r(R)$ whose elements $q$ satisfy $Jq\subseteq R$
for some dense left ideal $J$ of $R$ (see, e.g.~\cite{lam} for the definition of dense ideal). An alternative abstract
characterization of $Q_{\sigma}(R)$ is given in~\cite[Proposition
2.1]{lann}, as follows: Given a dense left ideal $I$ and a dense
right ideal $J$ of $R$, we say that a pair of maps $(f,g)$, where
$f\colon I\to R$ is a left $R$-homomorphism and $g\colon J\to R$
is a right $R$-homomorphism, is \emph{compatible} provided that
$f(x)y=xg(y)$ whenever $x\in I$ and $y\in J$.

Two sets of data $(f,g,I,J)$ and $(f',g',I',J')$ as above are
\emph{equivalent} if $f$ and $f'$ agree on some dense left ideal
contained in $I\cap I'$ and similarly for $g$ and $g'$. Denote by
$[f,g,I,J]$ the equivalence class of $(f,g,I,J)$ as above. Under
natural operations, the set of such equivalence classes is a ring
isomorphic to $Q_{\sigma}(R)$.

\begin{prop}
\label{20/02/05,4} Let $L\subseteq Q$ be a dense extension of Lie
algebras. Suppose that $Q$ is a multiplicatively semiprime algebra
of quotients of $L$. There is then an injective algebra
homomorphism $\tau\colon A(Q)\to Q_{\sigma}(A_0)$.
\end{prop}

\begin{proof}
Observe first that $L$ is (multiplicatively) semiprime by
Remark~\ref{20/02/05,3}. Next, given $\mu$ in
$A(Q)\setminus\{0\}$, we have proved in Theorem~\ref{elteorema}
that there is an ideal $I$ of $L$ with $\ann L(I)=0$ such that
$\mu \widetilde{I}+ \widetilde{I}\mu\subseteq A_0$ (where
$\widetilde I$ is the ideal of $A_Q(L)$ generated by the elements
of the form $\ad x$ for $x$ in $I$). By the observation made in
Remark~\ref{obs}, the left ideal $A_0\widetilde I+\widetilde I$
has zero right annihilator in $A_0$. By Corollary~\ref{20/02/05,6}
we have that $\lan{A (Q)}(\widetilde{I})=0$ and so the right ideal
$\widetilde IA_0+\widetilde I$ has zero left annihilator in $A_0$.
Since $A_0$ is semiprime by Lemma~\ref{20/02/05,2}, the ideals
$A_0\widetilde I+\widetilde I$ and $\widetilde IA_0+\widetilde I$
are left and right dense in $A_0$, respectively. Note that
$(A_0\widetilde I+\widetilde I)\mu\subseteq A_0$ and that
$\mu(\widetilde I A_0+\widetilde I)\subseteq A_0$. Hence, right
and left multiplication by $\mu$ produce two homomorphisms
$R_\mu\colon A_0\widetilde I+\widetilde I\to A_0$ and
$L_\mu\colon\widetilde I A_0+\widetilde I\to A_0$ which are left
(respectively, right) $A_0$-linear and compatible.

If $\mu=0$, we may then take $I=L$ and so $\widetilde I=A_Q(L)$,
which has zero right and left annihilator in $A(Q)$ (and hence in
$A_0$).

This allows us to define a map
\[
\tau\colon A(Q) \longrightarrow Q_\sigma(A_0)
\]
by $\tau(\mu)= [R_\mu, L_\mu, A_0
\widetilde{I}+\widetilde{I},\widetilde{I}A_0 +\widetilde{I}]$.
Note that, given any other pair of left and right dense ideals
$I'$ and $J'$ of $A_0$ such that $I'\mu\subseteq A_0$ and $\mu
J'\subseteq A_0$, we have that $\tau(\mu)=[R_\mu, L_\mu, I',J']$.
It follows easily from this that $\tau$ is an algebra
homomorphism.

In order to check that $\tau$ is injective, assume that
$\tau(\mu)=0$. This implies that there is a left dense ideal $I'$
of $A_0$ with zero right annihilator in $A_0$ such that $I'\mu=0$.
Now, $A_0$ is a left quotient algebra of $I'$ and by
Theorem~\ref{elteorema}, $A(Q)$ is a left quotient algebra of
$A_0$. Thus $A(Q)$ is a left quotient algebra of $I'$ and
therefore $I'\mu=0$ forces $\mu=0$.
\end{proof}

\begin{thm}
\label{20/02/05,7} Let $L\subseteq Q$ be a dense extension of Lie
algebras. Suppose that $Q$ is a multiplicatively semiprime algebra
of quotients of $L$. Then $A_0$ is semiprime and
$Q_\sigma(A_0)=Q_\sigma(A(Q))$.
\end{thm}

\begin{proof}
The semiprimeness of $A_0$ follows from Lemma~\ref{20/02/05,2}.
Now, apply Proposition~\ref{20/02/05,4} and~\cite[Theorem
2.5]{lann} to obtain $Q_\sigma(A_0)=Q_\sigma(A(Q))$.
\end{proof}

We close by exploring the possible converses to
Corollary~\ref{leftquotient} in the presence of dense extensions of Lie algebras.

\begin{definition} {\rm Given an extension of algebras
$A\subseteq S$, we say that $S$ is \emph{right ideally absorbed
into} $A$ if for any $q$ in $S\setminus\{0\}$ there is an ideal
$I$ of $A$ with  $\lan{A}(I)=0$ and such that $0\neq qI\subseteq
A$. \emph{Left ideally absorbed} can be defined analogously.}

\end{definition}

Observe that, in the definition above, we are requiring that the
ideal $I$ is two-sided rather than just a right ideal. In fact, in
the latter case this would be equivalent to saying that $S$ is a
right quotient algebra of $A$ (see~\cite[Lemma 2.14]{msm}).

\begin{prop}

\label{18/11/04,1} Let $L\subseteq Q$ be a dense extension of Lie
algebras with $Z(Q)=0$. Suppose that $A(Q)$ is \ria $A_0$. Then
$Q$ is an \qa of $L$.
\end{prop}

\begin{proof}
Let $q\in Q\setminus\{0\}$. Since $Z(Q)=0$, we have that $\ad
q\neq 0$. By hypothesis, there is an ideal $I$ of $A_0$ such that
$\lan{A_0}(I)=0$ and $0\neq\ad q I\subseteq A_0$. Set
$I_0=\{\alpha (x)\mid \alpha\in I\,, x\in L\}$, which is an ideal
of $L$ (see, e.g. the argument in the proof of
Theorem~\ref{elteorema}). Moreover, $\ann{L} (I_0)=0$. Indeed,
suppose that an element $x$ in $L$ satisfies $[x, I_0]=0$. By
definition, this means that $\ad x I(L)=0$, and since the
extension is dense we have that $\ad x I=0$. Thus $\ad x\in
\lan{A_0}(I)=0$. By Lemma~\ref{densext}, $Z(L)=0$ and so $x=0$.

Finally, $0\neq [q, I_0]=(\ad q I)(L)\subseteq A_0(L)\subseteq L$.
\end{proof}

\begin{cor} \label{18/11/04,2} Let $L\subseteq Q$ be a dense extension
of Lie algebras, with $L$ semiprime and $Z(Q)=0$. If $A(Q)$ is
\ria  $A_0$, then $A(Q)$ is \lia $A_Q(L)$.
\end{cor}

\begin{proof} By Proposition~\ref{18/11/04,1}, if $A(Q)$ is \ria $A_0$
we have that $Q$ is an \qa of $L$. We then may apply
Theorem~\ref{elteorema} to achieve the conclusion.
\end{proof}

\end{document}